# Polynomial Solutions of Differential Equations

H. Azad and M. T. Mustafa
Department of Mathematics and Statistics,
King Fahd University of Petroleum & Minerals, Dhahran, Saudi Arabia
hassanaz@kfupm.edu.sa, tmustafa@kfupm.edu.sa

**ABSTRACT.** We show that any differential operator of the form $L(y) = \sum_{k=0}^{k=N} a_k(x) y^{(k)}$, where $a_k$ is a real polynomial of degree $\leq k$, has all real eigenvalues in the space of polynomials of degree at most $n$, for all $n$. The eigenvalues are given by the coefficient of $x^n$ in $L(x^n)$.

If these eigenvalues are distinct, then there is a unique monic polynomial of degree $n$ which is an eigenfunction of the operator $L$- for every non-negative integer $n$. As an application we recover Bochner's classification of second order ODEs with polynomial coefficients and polynomial solutions, as well as a family of non-classical polynomials.

The subject of polynomial solutions of differential equations is a classical theme, going back to Routh [10] and Bochner [3]. A comprehensive survey of recent literature is given in [6]. One family of polynomials- namely the Romanovski polynomials [4, 9] is missing even in recent mathematics literature on the subject [8]; these polynomials are the main subject of some current Physics literature [9, 11]. Their existence and – under a mild condition - uniqueness and orthogonality follow from the following propositions. The proofs use elementary linear algebra and are suitable for class-room exposition. The same ideas work for higher order equations [1].

*Proposition 1*

*Let $L(y) = \sum_{k=0}^{k=N} a_k(x) y^{(k)}$, where $a_k$ is a real polynomial of degree $\leq k$. Then $L$ operates on the space $P_n$ of all polynomials of degree at most $n$. It has all real eigenvalues and the eigenvalues are given by the coefficient of $x^j$ in $L(x^j)$ for all $j \leq n$.*

*If the eigenvalues are distinct, then L has, up to a constant, a unique polynomial of every degree which is an eigenfunction of L*

*Proof:*

Let $L$ be as in the statement of the proposition. Since $L(x^n)$ is a sum of a multiple of $x^n$ plus lower order terms, it is clear that $L$ operates on every $P_j$, $j \leq n$. Therefore the eigenvalues are given by the coefficient of $x^j$ in $L(x^j)$ and $L$ has eigenfunctions in each $P_j$.

Assume that the eigenvalues of $L$ are distinct. Then $P_n$ has a basis of eigenfunctions and, for reasons of degree, there must be an eigenfunction of degree $n$, for every $n$. Therefore, up to a constant, there is a unique eigenfunction of degree $n$ for all $n$.

□

We now concentrate on second order operators, leaving the higher order case to [1]. Let $L(y) = a(x)y'' + b(x)y'$, where $\deg(a) \leq 2$, $\deg(b) \leq 1$. Following Bochner [3] if $\deg(a) = 2$ then by scaling and translation, we may assume that $a(x) = x^2 - 1, x^2 + 1$ or $x^2$. Applying the above proposition we then have the following result.

*Proposition 2*

(i) *The equation* $(x^2 + \varepsilon)y'' + (\alpha x + \beta)y' + \lambda y = 0$, $\varepsilon = 1, -1$ *has unique monic polynomial solutions in every degree if $\alpha > 0$ or if $\alpha < 0$ and it is not an integer. If $\alpha = -(n + k - 1)$ for $0 \leq k \leq (n-1)$, then the eigenspace in $P_n$ for eigenvalue $\lambda = n(n-1) + \alpha n$ is 2-dimensional.*

(ii) *The equation $xy'' + (\alpha x + \beta)y' + \lambda y = 0$ has unique monic polynomial solutions in every degree if $\alpha \neq 0$*

(iii) *The equation $y'' + (\alpha x + \beta)y' + \lambda y = 0$ has unique monic polynomial solutions in every degree if $\alpha \neq 0$*

In this proposition there is no claim to any kind of orthogonality properties. Nevertheless, the non-classical functions appearing here are of great interest in Physics and their properties and applications are investigated in [4, 9, 11].

The classical Legendre, Hermite, Laguerre and Jacobi make their appearance as soon as one searches for self-adjoint operators. Their existence and orthogonality properties [cf:8, p.80-106,2,7] can be obtained elegantly in the context of elementary Sturm-Liouville theory.

*Proposition 3*

*Let $L$ be the operator defined by $L(y) = a(x)y'' + b(x)y' + c(x)y$ on a linear space $C$ of functions which are at least two times differentiable on a finite interval $I$.*
*Define a bilinear function on $C$ by*

$$(y,u) = \int_I p(x) y(x) u(x) dx ,$$

*where $p$ is two times differentiable and non-negative and does not vanish identically in any subinterval of $I$.*
*Then*

$$(Ly,u)-(y,Lu) = pa(uy' - u'y) \Big|_\alpha^\beta$$

*if*

$$(pa)' = pb.$$

*Proof*:

Let $\alpha, \beta$ be the end points of $I$. So

$$(Ly,u) = \int_\alpha^\beta p(ay'' + by' + cy) u \, dx .$$

Using integration by parts, we find that $(Ly, u) - (y, Lu)$ will contain only boundary terms if $(pau)'' - (pbu)' = pau'' + pbu'$, for all $u$.

This simplifies to

$$[(pa)'' - (pb)']u + 2(pa)'u' = 2pbu'.$$

Equating coefficients of $u$ and $u'$ on both sides, we get the differential equations for $p$:

$$(pa)'' - (pb)' = 0 \text{ and } (pa)' = pb,$$

so in fact we need only the equation

$$(pa)' = pb.$$

The boundary terms now simplify to

$$(pau)y' - (pa)'uy - (pa)u'y + (pbu)y \Big|_\alpha^\beta = pa(uy' - u'y)\Big|_\alpha^\beta$$

The differential equation for the weight is $a'p + ap' = pb$, which integrates to

$$p = e^{\int \frac{(b-a')}{a}dx} = \frac{1}{|a|}e^{\int \frac{b}{a}dx}.$$

$\square$

*Examples*:

(1) **Jacobi polynomials**

First note that for any differentiable function $f$ with $f'$ continuous, the integral $\int_0^\varepsilon \frac{f(x)}{x^\alpha}dx$ is finite if $\alpha < 1$ - as one sees by using integration by parts.

Consider the equation $(1-x^2)y'' + (\alpha x + \beta)y' + \lambda y = 0$. As above, the weight function $p(x)$ for the operator

$$L(y) = (1-x^2)y'' + (ax+b)y'$$

is

$$p(x) = \frac{1}{1-x^2}e^{\int\left(\frac{\beta+\alpha}{2}\cdot\frac{1}{1-x} + \frac{\beta-\alpha}{2}\cdot\frac{1}{1+x}\right)dx} = \frac{(1+x)^{\frac{\beta-\alpha-2}{2}}}{(1-x)^{\frac{\beta+\alpha+2}{2}}} = \frac{1}{(1-x)^{\frac{\beta+\alpha+2}{2}}(1+x)^{\frac{-\beta+\alpha+2}{2}}}.$$

So $\int_{-1}^{1} p(x)f(x)dx$ would be finite if $\beta + \alpha < 0$ and $-\beta + \alpha < 0$, that is, if $\alpha < \beta < -\alpha$.

The weight is not differentiable at the end points of the interval. So, first consider $L$ operating on twice differentiable functions on the interval $[-1+\varepsilon, 1-\varepsilon]$. If $u, v$ are functions in this class then by *Proposition 3*,

$$\int_{-1+\varepsilon}^{1-\varepsilon} p(x)L(u(x))v(x)dx - \int_{-1+\varepsilon}^{1-\varepsilon} p(x)u(x)L(v(x))dx = p(x)a(x)(u(x)v'(x) - u'(x)v(x))\Big|_{-1+\varepsilon}^{1-\varepsilon}$$

Moreover, $(1-x^2)p(x) = (1-x)^{\frac{-(\beta+\alpha)}{2}}(1+x)^{\frac{\beta-\alpha}{2}}$ is continuous on the interval $[-1,1]$ and vanishes at the end-points -1 and 1. Therefore, if we define

$$(u,v) = \lim_{\varepsilon \to 0} \int_{-1+\varepsilon}^{1-\varepsilon} p(x)u(x)v(x)dx,$$ then $L$ would be a self-adjoint operator on all

polynomials of degree $n$ and so, there must be, up to a scalar, a unique polynomial which is an eigen function of $L$ for eigenvalue $-n(n-1) + n\alpha$.

So these polynomials satisfy the equation

$$(1-x^2)y'' + (\alpha x + \beta)y' + (n(n-1) - n\alpha)y = 0$$

and this equation has unique monic polynomial eigenfunctions of every degree, which are all orthogonal.

The Legendre and Chebyschev polynomials are special cases, corresponding to the values $\alpha = -1, -2, -3$ and $\beta = 0$.

(2) **The equation** $t(1-t)y'' + (1-t)y + \lambda y = 0$

This equation is investigated in [5] and the eigenvalues determined experimentally, by machine computations. Here, we will determine the eigenvalues in the framework provided by *Proposition* 3.

Let $L(y) = t(1-t)y'' + (1-t)y'$. Let $P_n$ be the space of al polynomials of degree at most $n$. As $L$ maps $P_n$ into itself, the eigenvalues of $L$ are given by the coefficient of $x^n$ in $L(x^n)$. The eigenvalues turn out to be $-n^2$. As these eigenvalues are distinct, there is, up to a constant, a unique polynomial of degree $n$ which is an eigenfunction of $L$.

The weight function is $p(t) = \dfrac{1}{|(1-t)|} = \dfrac{1}{(1-t)}$ on the interval $[0,1]$ and it is not integrable. However, as $L(y)(1) = 0$, the operator maps the space $V$ of all polynomials that are multiples of $(1-t)$ into itself. Moreover,

$\int_0^1 p(t)((1-t)\psi(t))^2 dt$ is finite.

The requirement for $L$ to be self-adjoint on $V$ is $t(\xi\eta' - \xi'\eta)\big|_0^1 = 0$ for all $\xi, \eta$ in $V$. As $\xi, \eta$ vanish at 1, the operator $L$ is indeed self-adjoint on $V$.

Let $V_n = (1-t)P_n$, where $P_n$ is the space of al polynomials of degree at most $n$.

As the codimension of $V_n$ in $V_{n+1}$ is 1, the operator $L$ must have an eigenvector in $V_n$ for all the degrees from 1 to $(n+1)$.

If $y = (1-t)\psi$ is an eigenfunction and $\deg(\psi) = n$ then, by the argument as in the examples above, we see that the corresponding eigenvalue is $\lambda = -(n+1)^2$.

Therefore, up to a scalar, there is a unique eigenfunction of degree $(n+1)$ which is a multiple of $(1-t)$ and all these functions are orthogonal for the weight $p(t) = \dfrac{1}{(1-t)}$. Using the uniqueness up to scalars of these functions, the eigenfunctions are determined by the differential equation and can be computed explicitly.

### (3) The Finite Orthogonality of Romanovski Polynomials

These polynomials are investigated in Refs [11,9] and their finite orthogonality is proved also proved there. Here, we establish this in the framework of *Proposition3*.

The Romanovski polynomials are eigenfunctions of the operator $L(y) = (1+x^2)y'' + (\alpha x + \beta)y'$. For $\alpha > 0$ or $\alpha < 0, \alpha$ not an integer, there is only one monic polynomial in every degree which is an eigenfunction of $L$; for $\alpha$ a non-positive integer, the eigenspaces can be 2 dimensional for certain degrees (*Propostion2*).

The formal weight function is $p(x) = (x^2+1)^{(\frac{\alpha-2}{2})} e^{\beta \tan^{-1}(x)} = (x^2+1)^{\frac{\gamma}{2}} e^{\beta \tan^{-1}(x)}$, where $\gamma = (\alpha - 2)$. Therefore, a polynomial of degree $N$ is integrable over the reals with weight $p$ if and only if $((N + \gamma + 1) < 0$ and if the product of two polynomials $P, Q$ is integrable, then the polynomials are themselves integrable for the weight $p$.

Arguing as in the proof of *Proposition* 3, we find that $(LP, Q) - (P, LQ) = (x^2+1)p(x)(PQ' - P'Q)\big|_{-\infty}^{\infty} = 0$, because the product

$(x^2 +1)p(x)(PQ' - P'Q)|$ is asymptotic to $x^{(2+\gamma+\deg(P)+\deg(Q)-1)} = x^{\deg(P)+\deg(Q)+\gamma+1}$ and $(\deg(P) + \deg(Q) + \gamma + 1) < 0$.

Therefore, if $P, Q$ are integrable eigenfunctions of $L$ with different eigenvalues and $(\deg(P) + \deg(Q) + \gamma + 1) < 0$, then $P, Q$ are orthogonal.

For several non-trivial applications to problems in Physics, the reader is referred to the paper [9].

***Conclusion***: In this note, which should have been written at least hundred years ago, we have rederived several results from classical and recent literature from a unified point of view by a straightforward application of basic linear algebra.

Some of these polynomials are not discussed in the standard textbooks on the subject, e.g. [8]- as pointed out in Ref [9].

We have also derived the orthogonality- classical as well as finite- of these polynomials from a unified point of view.